\documentclass[12pt]{article}
\usepackage{verbatim}
\newcommand{\nocopyright}{
No Copyright\thanks{
The authors hereby waive all copyright
and related or neighboring rights to this work,
and dedicate it to the public domain.
This applies worldwide.
}}
\author{
Peter G. Doyle \and Jean Steiner
}
\title{Commuting time geometry of ergodic Markov
chains}
\date{Version 2A3 dated 10 October 2017
\\ \nocopyright
}
\usepackage{tensor}
\newcommand{\mod}{/}

\newcommand{\Prob}{\mathrm{Prob}}

\newcommand{\Lap}{\Delta}

\newcommand{\inv}{{-1}}
\newcommand{\invstar}{{-1 \star}}
\newcommand{\qed}{\rule{2mm}{2.5mm}}

\newcommand{\phibar}{\bar{\phi}}
\newcommand{\psibar}{\bar{\psi}}
\newcommand{\one}{\mathrm{1}}

\newcommand{\Pinf}{{P^{(\infty)}}}
\newcommand{\tdelta}{\tensor{\delta}}
\newcommand{\half}{\frac{1}{2}}
\newcommand{\bra}{\langle}
\newcommand{\ket}{\rangle}
\newcommand{\phiinv}{{\phi^{-1}}}
\newcommand{\phistar}{{\phi^\star}}
\newcommand{\phiinvstar}{{\phi^{-1 \star}}}
\newcommand{\Vstar}{{V^\star}}
\newcommand{\Vstarstar}{{V^{\star \star}}}
\newcommand{\R}{{\bf R}}
\newcommand{\after}{\circ}

\begin{document}
\maketitle

\begin{abstract}

We show how to map the states of an ergodic Markov chain to Euclidean
space so that the squared distance between states is the expected
commuting time,
and we find a minimax characterization of commuting times.
These results are familiar in the case of time-reversible
chains, where techniques of classical electrical theory apply,
and the minimax simplifies to a minimum principle.
In presenting these results, we take the opportunity
to develop Markov chain theory in a `conformally correct' way.

\end{abstract}

\section{Overview}

In an eye-opening paper,
Chandra, Raghavan, Ruzzo, Smolensky, and Tiwari
\cite{crrst:commute}
revealed the central importance of expected commuting times
for the theory of time-reversible Markov chains.
Here we extend the discussion to general, non-time-reversible chains.

We begin by showing how to embed the states in a Euclidean space
so that the squared distance between states is the commuting time.
In the time-reversible case,
Leibon et al.
have used Euclidean embeddings to great effect
as a way to visualize a chain,
and reveal natural clustering of states.
Our embedding theorem shows that non-time-reversible chains should
be amenable to the same treatment.

Looking beyond the Euclidean embedding,
we find a natural minimax characterization of commuting times.
For time-reversible chains, this simplifies
to the minimum dissipation principle familiar from electrical theory.

In presenting these results, we will be taking a `conformally correct'
approach to Markov chains.
Briefly, a conformal change to a Markov chain changes its equilibrium
measure, but not its equilibrium transition rates.
The opportunity to develop this conformally correct approach is at
least as important to us as the particular results
we'll be discussing here.

\subsubsection*{CORRECTION}
In the first version of this paper, we (Peter) claimed that
the minimax characterization of commuting times
implies the following monotonicity law:
If all equilibrium interstate transition rates are increased,
then all commuting times are diminished.
For time-reversible chains, this monotonicity law is an ancient
and powerful tool.
As Russ Lyons and Yuval Perez have kindly pointed out,
the generalization to non-reversible chains fails miserably
already for a merry-go-round chain, where 
from $i\, (\mathrm{mod}\ n)$ you move to $i+1\,(\mathrm{mod}\ n)$
with probability $1/2$,
and otherwise stay where you are.
The commuting time between any two distinct states is $2n$,
while for a simple random walk on the integers mod $n$
all equilibrium transition rates are larger,
but commuting times between distant states are quadratic in $n$.

The good news is that we had found no applications for this bogus
generalization of monotonicity,
and had even observed,
`it is questionable how useful it will prove in the general case.'
Cold comfort, but comfort nonetheless.

Russ and Yuval also settled the question we had left open about whether
any commuting time matrix can arise from a reversible chain:
They showed that the answer is no.
We'll come back to this.

\section{The problem}

The \emph{commuting time} $T_{ab}$
between two states $a,b$ of an ergodic Markov chain
is the expected time, starting from $a$, to go to $b$ and then back to $a$.
Evidently $T_{ab}=T_{ba}$ and
\[
T_{ac} \leq T_{ab}+T_{bc}
.
\]
Thus it might seem natural to think of $T_{ab}$ as a measure of the distance
between $a$ and $b$.
But in fact it is most natural to think of $T_{ab}$ as the
\emph{squared distance} between $a$ and $b$.
The reason is that, as we will see,
there is a natural way to identify the states of the chain
with points in a Euclidean
space having quadratic form $||x||^2$ such that for any states $a,b$
we have
\[
T_{ab} = ||a-b||^2
.
\]
Now that we are interpreting $T_{ab}$
as a squared distance,
the inequality
$T_{ac} \leq T_{ab}+T_{bc}$
tells us that
\[
||a-c||^2 \leq ||a-b||^2 +||b-c||^2
.
\]
This means that all angles $\angle abc$ are acute
(at least weakly: some might be right angles).

Realizing commuting times as squared distances
is straight-forward for time-reversible chains.
Here's a sketch, meant only for orientation:
We won't rely on any of this below.
Time-reversible chains correspond exactly to resistor networks,
with $T_{ab}$ corresponding to the effective resistance
between $a$ and $b$.
This effective resistance is the energy of a unit current flow
from $a$ to $b$.
The energy of a flow is its squared distance with respect to the
energy norm on flows.
If we associate to state $i$ the unit current
flow from $i$ to some arbitrary reference
vertex (the `ground'),
then the difference between the flows associated to $a$ and $b$ will
be the unit current flow from $a$ to $b$,
having square norm $T_{ab}$.

The trick will be to extend this result to non-time-reversible chains.

{\bf Note.}
If it were the case
that to any chain there corresponds a time-reversible
chain having the same $T$, up to multiplication by a positive constant,
this extension would be immediate.
It is easy enough to compute what the transition rates of this
time-reversible chain would have to be,
the question is whether these are necessarily positive.
Russ Lyons found by computer search that the answer is no:
There are chains whose commuting-time matrix
does not arise from a time-reversible chain.
In fact such chains become plentiful as the number of states increases.

Before proceeding, we should observe
that the triangle inequality for squared lengths
is not in itself
a sufficient condition for realizability of a Euclidean simplex.
It \emph{is} sufficient for tetrahedra (four vertices in 3-space),
but for five vertices we have the following
counterexample.
Take
\[
T=
\left(
\begin{array}{ccccc}
0& 7& 7& 7& 13
\\
7& 0& 12& 12& 7
\\
7& 12& 0& 12& 7
\\
7& 12& 12& 0& 7
\\
13& 7& 7& 7& 0
\end{array}
\right)
\]
This matrix
is not realizable because the associated quadratic form with matrix
\[
\half
\left(
\begin{array}{llll}
 14 & 2 & 2 & 13 \\
 2 & 14 & 2 & 13 \\
 2 & 2 & 14 & 13 \\
 13 & 13 & 13 & 26
\end{array}
\right)
\]
is not positive definite:
It has the eigenvalue
$\half (22-\sqrt{523}) \approx -0.434597$.
Since we're going to see that commuting time matrices are always realizable,
this means in particular that this matrix $T$
cannot arise as the matrix of commuting times of a Markov chain.

\section{The short answer}

Below we will give the honest solution to this problem,
developing in a thoroughgoing way
what we will call the `conformally correct'
approach to Markov chains.
Here we just extract the answer to our embedding question,
and present it in a way that should be immediately accessible to
those familiar with
the standard theory of Markov chains,
as developed for example in
Grinstead and Snell
\cite{grinsteadSnell:prob}.
The only caveat is that we will be using tensor notation, i.e. writing
some indices up rather than down.
You can look at section \ref{sec:tensor} below for remarks about this,
but if you prefer you can just view this as an idiosyncracy,
as long as you bear in mind that
$\tensor{Z}{_i^j}$ represents a different array of numbers from $Z_{ij}$.

Consider a discrete-time Markov chain with transition probabilities
\[
\tensor{P}{_i^j} = \Prob(\mbox{next at $j$}|\mbox{start at $i$})
.
\]
Assume the chain is ergodic so there is a unique equilibrium measure $w^i$
with
\[
\sum_i w^i \tensor{P}{_i^j} = w^j
\]
and
\[
\sum_i w^i = 1
.
\]
Define
\[
\Lap^{ij} = w^i(\tensor{I}{_i^j}-\tensor{P}{_i^j})
,
\]
and note that
\[
\sum_i \Lap^{ij} = \sum_j \Lap^{ij} = 0
.
\]

Now define
\[
\tensor{Z}{_i^j}
=
(\tensor{I}{_i^j}-w^j)
+
(\tensor{P}{_i^j}-w^j)
+
(\tensor{{P^{(2)}}}{_i^j}-w^j)
+
\ldots
,
\]
where
$\tensor{{P^{(2)}}}{_i^j}
=
\sum_k \tensor{P}{_i^k} \tensor{P}{_k^j}$
represents the matrix square of $\tensor{P}{_i^j}$,
and the elided terms involve higher matrix powers.
Using conventional matrix notation
if we define $\tensor{\Pinf}{_i^j} = w^j$ we can write
\begin{eqnarray*}
Z
&=&
(I-\Pinf) + (P-\Pinf) + (P^{(2)}-\Pinf) + \ldots
\\&=&
(I-P+\Pinf)^{-1} - \Pinf
.
\end{eqnarray*}
(Note that Grinstead and Snell
\cite{grinsteadSnell:prob}
use the alternate definition
$Z=(I-P+\Pinf)^{-1}$,
which is less congenial
but works just as well in this context.)

Set
\[
Z_{ij} = \frac{1}{w^j}\tensor{Z}{_i^j}
.
\]
$Z_{ij}$ acts like an inverse to $\Lap^{ij}$ in the sense that
for any $u^i$ with $\sum_i u^i=0$,
we have
\[
\sum_{jk} u^j Z_{jk} \Lap^{kl} = u^l
\]
and
\[
\sum_{jk} \Lap^{ij} Z_{jk} u^k = u^i
.
\]

Standard Markov chain theory tells us that the expected time $M_{ab}$
to hit state $b$ starting from state $a$ is
\[
M_{ab} = Z_{bb} - Z_{ab}
.
\]
So for the commuting time we have
\[
T_{ab}=M_{ab}+M_{ba}
=
Z_{aa}-Z_{ab}-Z_{ba}-Z_{bb}
.
\]

For a vector $x=(x_i)_{i=1,\ldots,n}$
define
\[
||x||^2 =
\sum_{ij} x_i \Lap^{ij} x_j
.
\]
Please note that this does not make $\Lap^{ij}$ the matrix
of the quadratic form in the usual sense,
because in general
$\Lap^{ij} \neq \Lap^{ji}$.
The matrix of the form in the usual sense is the symmetrized version
$\half(\Lap^{ij}+\Lap^{ji})$.

Because
\[
\sum_i \Lap^{ij} = \sum_j \Lap^{ij} = 0
\]
we have the key identity
\[
||x||^2 = - \half \sum_{ij} \Lap^{ij} (x_i-x_j)^2
.
\]
Recalling the definition of $\Lap^{ij}$ gives
\[
||x||^2 =
\half \sum_{ij} w^i \tensor{P}{_i^j} (x_i-x_j)^2
.
\]
Thus the quadratic form $||x||^2$ is weakly positive definite,
but not strictly so, because it vanishes for constant vectors:
\[
||(c,\ldots,c)||^2 = 0
.
\]
It becomes strictly positive definite
if we identify vectors differing by a constant vector:
\[
(x_i)_{i=1,\ldots,n} \equiv (z_i+c)_{i=1,\ldots,n}
.
\]
This Euclidean space (vectors mod constant vectors, with the
pushed-down quadratic form) is where we will embed our chain.

To get the embedding, map state $a$ to the vector
\[
f(a) = (Z_{ai})_{i=1,\ldots,n}
.
\]
For the difference between the images of $a$ and $b$ we have
\[
(f(a)-f(b))_i = Z_{ai}-Z_{bi}
=
\sum_k 
(\tdelta{_a^k}-\tdelta{_b^k}) Z_{ki}
,
\]
with $\tdelta{_i^j}$ the Kronecker delta.
We want to see that $f(a)-f(b)$ has square norm $T_{ab}$.

From the generalized inverse relationship between $Z_{ij}$ and $\Lap^{ij}$
and the fact that
\[
\sum_k
(\tdelta{_a^k}-\tdelta{_b^k}) = 0
\]
we have
\[
\sum_{i} (Z_{ai}-Z_{bi}) \Lap^{ij}
=
\sum_{ki} (\tdelta{_a^k}-\tdelta{_b^k}) Z_{ki} \Lap^{ij}
=
\tdelta{_a^j}-\tdelta{_b^j}
.\]
So
\begin{eqnarray*}
||f(a)-f(b)||^2
&=&
\sum_{ij} (Z_{ai}-Z_{bi}) \Lap^{ij} (Z_{aj}-Z_{bj})
\\&=&
\sum_j (\tdelta{_a^j}-\tdelta{_b^j}) (Z_{aj}-Z_{bj})
\\&=&
Z_{aa}-Z_{ab}-Z_{ba}+Z_{bb}
\\&=&
T_{ab}
.
\end{eqnarray*}
There you have it.

\section{What just happened}

We want to explain the proof we have just given in more conceptual terms.

Let $V$ be a finite-dimensional real vector space,
and $\Vstar$ the dual space,
consisting of linear functionals $\phi:V \to \R$.
For $u \in \Vstar$, $x \in V$ write
\[
\bra u,x \ket_V = u(x)
\]
for the natural pairing between $V$ and $\Vstar$.
Identify $V$ with $\Vstarstar$ as usual:
\[
\bra x,u \ket_\Vstar
=
u(x)
=
\bra u,x \ket_V
.
\]
To a map $f:V \to W$ we associate the adjoint map
$f^\star:W^\star \to V^\star$, such that for $u \in W^\star$, $x \in V$
\[
\bra f^\star(u),x \ket_V = \bra u, f(x) \ket_W = u(f(x))
.
\]

A bilinear form on $V$ arises from a linear map
\[
\phi:V \to \Vstar
\]
via
\[
L_\phi(x,y) = \bra \phi(x),y \ket_V
.
\]
The adjoint map
\[
\phistar: V \to \Vstar
\]
yields the transposed bilinear form
\[
L_\phistar(x,y) 
=
\bra \phistar(x),y \ket_V
=
\bra x,\phi(y) \ket_\Vstar
=
\bra \phi(y),x \ket_V
=
L_\phi(y,x)
.
\]

If $\phi$ is invertible the inverse
\[
\phiinv: \Vstar \to V
\]
yields the form $L_\phiinv$ on $\Vstar$:
\[
L_\phiinv(u,v)
=
\bra \phi^{-1}(u),v \ket_\Vstar
=
\bra v, \phiinv(u) \ket_V
.
\]

The forms $L_\phistar$ and $L_\phiinv$ are conjugate, because
\[
L_\phiinv(u,v)
=
\bra v,\phiinv(u) \ket_V
=
L_\phi(\phiinv(v),\phiinv(u))
=
L_\phistar(\phiinv(u),\phiinv(v))
.
\]
Going back the other way,
\[
L_\phistar(x,y)
=
L_\phiinv(\phi(x),\phi(y))
.
\]

From these two equations, we get two distinct ways to conjugate
$L_\phi$ to $L_\phiinvstar$.
Plugging $\phi=(\phiinv)^{-1}$ into the first
and putting $(x,y)$ for $(u,v)$, we get
\[
L_\phi(x,y)=L_\phiinvstar(\phi(x),\phi(y))
.
\]
Plugging $\phi=(\phistar)^\star$ into the second we get
\[
L_\phi(x,y)=L_\phiinvstar(\phistar(x),\phistar(y))
.
\]
Now putting $\phistar$ for $\phi$ we see that 
in fact there were two ways to conjugate $L_\phiinv$
to $L_\phistar$:
\[
L_\phistar(x,y)
=
L_\phiinv(\phi(x),\phi(y))
=
L_\phiinv(\phistar(x),\phistar(y))
.
\]

Having two ways to conjugate $L_\phi$ to $L_\phiinvstar$ gives us
an automorphism $\phiinv \after \phistar$ of $L_\phi$:
\[
L_\phi(x,y)=
L_\phi(\phiinv(\phistar(x)),\phiinv(\phistar(y)))
.
\]
Along with $\phiinv \after \phistar$ we also have the inverse
automorphism $\phiinvstar \after \phi$:
\[
L_\phi(x,y)=
L_\phi(\phiinvstar(\phi(x)),\phiinvstar(\phi(y)))
.
\]
We could also consider powers other than $-1$ of our automorphism,
but we don't need to, because the conjugacy between $L_\phi$ and
$L_\phistar$ is canonical (in the sense of being equivariant with
respect to taking duals and inverses) up to this factor of two.
The difference between them,
as measured by the automorphism $\phiinv \after \phistar$,
measures the antisymmetry of $L_\phi$.
It is destined to play an important role in our future.

Looking now at the level of quadratic forms $Q_\phi(x)=L_\phi(x,x)$,
everything in sight is conjugate:
\[
Q_\phi(x)=Q_\phistar(x)
;
\]
\[
Q_\phiinv(u) = Q_\phiinvstar(u)
= Q_\phi(\phiinv(u)) = Q_\phi(\phiinvstar(u))
.
\]

All this nonsense can be made much more concrete using matrices.
Let $V = \R^n$ and represent $x \in V$, $u \in \Vstar$ as column and
row vectors respectively, so that the pairing is just multplying
a row vector by a column vector:
\[
\bra u, x \ket_V = ux
.
\]
Denote transposition of matrices by $\star$.
Write
\[
L_\phi(x,y) = x^\star A y
,
\]
so that
\[
\phi(x) = x^\star A = (A^\star x)^\star
.
\]
Now
\[
\phiinv(u) = (u A^{-1})^\star = A^\invstar u^\star
,
\]
so
\[
L_\phiinv(u,v) =
\bra v, \phiinv(u) \ket_V
= v A^\invstar u^\star
= u A^{-1} v^\star
.
\]
Good!

Now to see the two conjugacies of $L_\phistar$ with $L_\phiinv$:
\[
A^\star A^\inv A = A^\star
;
\]
\[
A A^\inv A^\star = A^\star
.
\]
These combine to give two automorphisms of $L_\phi$:
\[
(A^\inv A^\star)^\star A (A^\inv A^\star)
=
A A^\invstar A A^\inv A^\star
=
A
;
\]
\[
(A^\invstar A)^\star A (A^\invstar A)
=
A^\star A^\inv A^\invstar A
=
A
.
\]
Hmm.  Why didn't we do it this way in the first place?

So, here's what happened with our Markov chain.
We started with the space $V=\R^n \mod \one$ with quadratic form
$L_\phi(x,y) = \sum_{ij} x_i \Lap^{ij} y_j$,
embedded the states in $\Vstar = \R^n \perp \one$ with quadratic form
$L_\phiinv(u,v) = \sum_{ij} u^i Z_{ij} v^j$,
and proved that $L_\phiinv$ is positive definite by showing that it
is conjugate to $L_\phi$.

\section{Tensor notation for Markov chains}\label{sec:tensor}

As you will already have noticed,
we are using tensor notation,
rather than trying to work within the confines of matrix notation,
as is usual in the theory of Markov chains.
For our purposes,
a tensor may be viewed as an array where 
some of the indices are written
as superscripts rather than subscripts.
Thus, for example, we write the transition rates for
a Markov chain as $\tensor{P}{_i^j}$,
and the equilbrium measure as $w^i$.

Where the indices of a tensor are placed makes a difference:
Thus $\tensor{Z}{_i^j}$ represents a different array from
$Z_{ij}$.
We may `raise' and `lower' these indices
as is usual with tensors,
though in this case the procedure is simpler than usual,
because to raise or lower an index $i$ we just multiply or divide by
the entries of $w^i$.
Thus we get $Z_{ij}$ from  $\tensor{Z}{_i^j}$ by lowering the index $j$:
\[
Z_{ij} = \frac{1}{w^i} \tensor{Z}{_i^j}
.
\]
We get back to $\tensor{Z}{_i^j}$ from $Z_{ij}$ by raising the index $j$:
\[
\tensor{Z}{_i^j}
=
w^j Z_{ij}
.
\]

We will still be able to use matrix notation to multiply
matrices (two-index tensors) and vectors (one-index tensors).
The beautiful thing is that when we do this, the indices take care of
themselves,
as long as the indices that get summed over when multiplying matrices
are paired high with low.
To show by example what this means,
if we write $C=AB$, it will entail (among other things) that
\[
\tensor{C}{_i^j} = 
\tensor{(AB)}{_i^j} =
\sum_k \tensor{A}{_i^k}\tensor{B}{_k^j}
=
\sum_k \tensor{A}{_i_k}\tensor{B}{^k^j}
,
\]
and 
\[
\tensor{C}{_i_j} = 
\tensor{(AB)}{_i_j} =
\sum_k \tensor{A}{_i^k}\tensor{B}{_k_j}
=
\sum_k \tensor{A}{_i_k}\tensor{B}{^k_j}
=
\sum_k \tensor{A}{_i_k} w^k \tensor{B}{_k_j}
=
\sum_k \tensor{A}{_i_k}\tensor{B}{^k^j}\frac{1}{w^j}
.
\]

{\bf Note.}
If you're familiiar with the Einstein summation convention,
be aware that we don't use it here.
It wouldn't work well in this context,
because we want to write $w^i Z_{ij}$
without automatically summing over $i$.
Fortunately,
for our purposes, using the notation of matrix multiplication turns out
to be even more convenient than the summation convention.

\section{What it means to be conformally correct}

We have said that we want our approach to be `conformally correct'.
Before we go further, a word about what this means.
(Skip this if you don't care.)

Conformal equivalence of Markov chains is most natural for continuous
time chains.
In that context two chains with transition rates $\tensor{A}{_i^j}$
and $\tensor{B}{_i^j}$
are conformally equivalent if
\[
\tensor{B}{_i^j} = \frac{1}{a_i} \tensor{A}{_i^j}
\]
where all $a_i>0$.
Generally we will also want the additional condition that
$\sum_i w^i a_i=1$ where $w^i$ is the equilibrium probability of being at $i$
for the $A$ chain.
With this `volume condition' the equilibrium probability of being at $i$
for the $B$ chain will be $w^i a_i$ and
\[
B^{ij} = w^i a_i \tensor{B}{_i^j} = w^i a_i \frac{1}{a_i} \tensor{A}{_i^j}
= A^{ij}
.
\]
Thus while the raw transition rates $\tensor{A}{_i^j}$ are not
conformal invariants, when we raise the index $i$ we get a new array
$A^{ij} = w^i \tensor{A}{_i^j}$ whose entries are conformal invariants:
They tell the rate at which transitions are made from $i$ to $j$ when
the chain is in equilibrium.

It is possible to talk about conformal equivalence of discrete time chains,
but it is not as pleasant as for continuous-time chains.
This is true so often in the theory of Markov chains!
And yet, for simplicity, we want to talk about discrete-time chains.
So our approach will be to do everything in such a way that the discussion
would be conformally invariant when translated from discrete to continuous
time.

So that's what it means for chains to be conformally equivalent.
As for `conformal correctness', we mean an approach that seeks to identify
and emphasize quantities that are conformally invariant.
And why should we do this?
Because it will pay.

\section{Visualizing commuting times}

One way to determine the expected commuting time $T_{ab}$ between $a$ and $b$
is to run the chain for a long time $T$
(beware of confusion!), paying attention to when the chain
is at $a$ or $b$ and ignoring other states.
If $R$ is the number of runs of $a$'s (which is within $1$ of the number of
runs of $b$'s),
then
\[
T_{ab} \approx T/R
.
\]
To keep track of $R$ we imagine painting our Markovian particle green when
it reaches $a$ and red when it reaches $b$.
Let $r_{ab}$ be the equilibrium rate at which red particles are being painted
green.
Ignoring end effects,
over our long time interval $T$,
$R$ above is the number of times a red particle gets painted green, thus
roughly $T r_{ab}$,
and it follows that
\[
T_{ab} = \frac{1}{r_{ab}}
.
\]
This is an instance of the general principle
from renewal theory
that when events happen at rate $r$,
the expected time between events is $1/r$.

{\bf Note.}  This painting business is very close to a model
developed by Kingman
\cite{kingman:paint}
and Kelly
\cite{kelly:paint}.
(See exercise 1 in section 3.3 of Doyle and Snell
\cite{doylesnell:walks}.)
However, I don't know that
Kingman and Kelley ever made the connection to commuting times,
and it is possible that their discussion concerned only time-reversible
chains.
Somebody should check this.

It is high time to observe that if $\hat{T}_{ab}$ is the
commuting time for the time-reversed chain
(according to the general convention that time-reversed quantities wear hats),
we have
\[
T_{ab}=T_{ba}=\hat{T}_{ab}=\hat{T}_{ba}
.
\]
We claim to be able to see this from our way of approximating $T_{ab}$
by observing the chain over a long time.
If we reverse a record of the chain moving forward for a long time,
we see roughly a record of the time-reversed chain starting in equlibrium.
In fact if we started the original chain in equilibrium we're golden.
If we started the chain not in equilibirum
(e.g. by starting at $a$,
as we might well be tempted to do),
there will be problems toward the
end of the time-reversed record,
as the time-reversed chain gets drawn to end where the forward chain began.
But this effect is negligible when $T$ is large.

\section{The Laplacian and the cross-potential}

Consider a discrete-time Markov chain with transition probabilities
\[
\tensor{P}{_i^j} = \Prob(\mbox{next at $j$}|\mbox{start at $i$})
.
\]
Assume the chain is ergodic,
so that there is a unique equilibrium measure $w^i$
with
\[
\sum_i w^i \tensor{P}{_i^j} = w^j
,
\]
\[
\sum_i w^i = 1
.
\]

Define the \emph{Laplacian}
\[
\Lap^{ij} = w^i(\tensor{I}{_i^j}-\tensor{P}{_i^j})
.
\]
For $i \neq j$, $-\Lap^{ij}$ tells
the equilibrium rate of transitions from $i$ to $j$;
$\Lap^{ii}$ tells the total rate of transitions
to and from states other than $i$.
The time-reversed Markov chain has Laplacian
$\hat{\Lap}^{ij} = \Lap^{ji}$.
A time-reversible chain has
$\Lap^{ij}=\Lap^{ji}$.

We have
\[
\sum_i \Lap^{ij} = \sum_j \Lap^{ij} = 0
.
\]
So considered as a matrix, $\Lap^{ij}$ is not invertible.
However, it has a generalized inverse $Z_{ij}$ with the property that
for any measure of total mass 0, which is to say
for any $u^i$ with $\sum_i u^i=0$,
we have
\[
\sum_{jk} u^j Z_{jk} \Lap^{kl} = u^l
\]
and
\[
\sum_{jk} \Lap^{ij} Z_{jk} u^k = u^i
.
\]
An equivalent way to write this is
\[
\sum_{jk} \Lap^{ij} Z_{jk} \Lap^{kl} = \Lap^{il}
,
\]
because if we think of $\Lap^{ij}$ as a matrix, its rows and columns both
span the space of measures with total mass 0.

A sensible choice for the generalized inverse $Z_{ij}$ is
\[
Z_{ij} = \frac{1}{w^j} \tensor{Z}{_i^j}
\]
where
\[
\tensor{Z}{_i^j}
=
(\tensor{I}{_i^j}-w^j)
+
(\tensor{P}{_i^j}-w^j)
+
(\tensor{{P^{(2)}}}{_i^j}-w^j)
+
\ldots
,
\]
where
$\tensor{{P^{(2)}}}{_i^j}
=
\sum_k \tensor{P}{_i^k} \tensor{P}{_k^j}$
represents the matrix square of $\tensor{P}{_i^j}$,
and the elided terms involve higher matrix powers.
Define $\tensor{\Pinf}{_i^j} = w^j$,
to suggest that the `infinitieth power' of $\tensor{P}{_i^j}$ has
all rows equal to the vector $w^i$.
We can write
\begin{eqnarray*}
Z
&=&
(I-\Pinf) + (P-\Pinf) + (P^{(2)}-\Pinf) + \ldots
\\&=&
(I-P+\Pinf)^{-1} - \Pinf
.
\end{eqnarray*}
This naturally translates into the formula we've given for
$\tensor{Z}{_i^j}$,
and from there, by `lowering the index j',
we get $Z_{ij}$.

For this choice of $Z$ we have the natural interpretation that
$\tensor{Z}{_i^j}$
is the expected excess number of visits to $j$ for a chain starting
at $i$ compared to a chain starting in equilibrium.
For the time-reversed chain we get
\[
\tensor{\hat{Z}}{_{ij}}=\tensor{Z}{_{ji}}
,
\]
and so in particular if the chain is time-reversible we have
$Z_{ij}=Z_{ji}$.

This is all very well, but we still do not want to prescribe this
particular choice of $Z$ because it is not conformally invariant:
It depends on the equilibrium measure $w^i$, and not just on the
Laplacian `matrix' $\Lap^{ij}$.
This makes it insufficiently canonical for us.

What \emph{is} canonical is the bilinear form
\[
B(u,v)= \sum_{ij} u^i Z_{ij} v^j
\]
when $u$ and $v$ are restricted to the subspace $S$ of measures of total
mass $0$:
\[
S = \{u^i: \sum_i u^i=0 \}
\]
Fixing $a,b,c,d$ and setting 
\[
u=\tdelta{_a^i}-\tdelta{_b^i}
;\;
v=\tdelta{_c^i}-\tdelta{_d^i}
\]
gives us the \emph{cross-potential}
\[
N_{abcd}
=
B(\tdelta{_a^i}-\tdelta{_b^i},\tdelta{_c^i}-\tdelta{_d^i})
=
Z_{ac}-Z_{ad}-Z_{bc}+Z_{bd}
.
\]
$N$ satisfies
\[
N_{bacd}=N_{abdc}=-N_{abcd}
.
\]
For the time-reversed process
\[
\hat{N}_{abcd}=N_{cdab}
.
\]

Clearly, knowing $N$ is the same as knowing $B$, or $\Lap$.
If we know $w$ as well as $N$ we can recover our
sensible-but-not-canonical $Z$:
\[
Z_{ij}=\sum_{kl} N_{ikjl}w^k w^l
.
\]
Different choices of $w$ in this formula lead to different $Z$'s,
but they all determine the same bilinear form $B$.
From $Z$ and $w$ we can recover $P$.

In general, it is useful to think of an ergodic Markov chain
as specified by the cross-potential $N$, which determines its
conformally invariant properties,
together with the equilibrium measure $w$.
Expressing formulas in these terms allows us to see
the extent to which quantities are conformally invariant
(like $N$, $B$, and $\Lap$) or not (like $w$, $Z$, $P$).

{\bf Complaint.}
$N$ and $w$ together don't quite determine the original transition rates
for a continuous-time Markov chain,
or rather, they wouldn't do so if we had some way to distinguish between
remaining at $i$ and moving from $i$ to $i$.
Such a distinction is not possible for discrete-time chains
represented by matrices,
but we could handle it in the continuous case by allowing for
non-zero transition rates on the diagonal.
Better yet, we could
reformulate Markov chain theory in the context of queuing networks
based on $1$-complexes
(graphs where loops and multiple edges are allowed).
This would give us a way to distinguish different ways of stepping
from $i$ to $j$.
A further step would be to allow a general distribution for
the time it takes to make a transition for $i$ to $j$.
This would be very helpful
when watching the chain only when it is in a subset
of its states,
as in the case above where we contemplated watching the chain only
when it is at $a$ or $b$.
We didn't say just what we meant by this, because it doesn't conveniently
fit into the usual formulation of Markov chain theory.

\section{Probabilistic and electrical interpretation}

We may interpret $N_{abcd}$ probabilistically as 
the equilibrium concentration difference between
$c$ and $d$ due to a unit flow of particles entering at $a$ and leaving at $b$.
Here's what this means.
Introduce Markovian particles at $a$ at a unit rate,
and remove them when they reach $b$.
Write the `dynamic equilibrium' measure of particles at $i$ as $w^i \phi_i$, so
that $\phi_i$ tells the concentration of particles relative
to the `static equilibrium' measure $w^i$.
Conservation of particles implies that
\[
w^i \phi_i \sum_j \tensor{P}{_i^j}
- \sum_j w^j \phi_j \tensor{P}{_j^i}
= \tdelta{_a^i}-\tdelta{_b^i}
.
\]
We hasten to rewrite this in the conformally correct form
\[
\sum_j \phi_j \Lap^{ji} = \tdelta{_a^i}-\tdelta{_b^i}
.
\]
Since also
\[
\sum_j (Z_{aj}-Z_{bj}) \Lap^{ji} = \tdelta{_a^i}-\tdelta{_b^i}
\]
and since the Laplacian $\Lap$ kills only constants, if follows that
\[
\phi_j = Z_{aj}-Z_{bj} + C
,
\]
and thus
\[
\phi_c - \phi_d = Z_{ac}-Z_{bc}-Z_{ad}+Z_{bd} = N_{abcd}
.
\]

From this probabilistic interpretation of $N$
we can see that $N_{abab} = C_{ab}$,
the commuting time between $a$ and $b$.
Indeed, 
in the particle-painting scenario introduced earlier,
$C_{ab}$ is the reciprocal of the
rate at which red particles are turning green at $a$.
Paying attention only to green particles, we
see green particles appearing at $a$ at rate $1/C_{ab}$,
and disappearing at $b$.
The equilibrium concentration of green particles at $i$ is the probability
$p_i$ of hitting $a$ before $b$ for the time-reversed chain,
and in particular $p_a=1$ and $p_b=0$, so the concentration difference
between $a$ and $b$ is $1$.
Multiplying the green flow by $C_{ab}$ normalizes it to a unit flow
with concentration difference $C_{ab}$ between $a$ and $b$.
So 
\[
C_{ab}=N_{abab}
.
\]

If we embellish this probabilistic scenario by
imagining that our particles carry a positive charge,
we may identify the net flow of particles with electrical current;
the concentration of particles 
(relative to the equilibrium measure) with electrical potential;
and differences of concentration with voltage drop.
With this terminology,
$N_{abcd}$ 
tells the voltage drop between $c$ and $d$ due to a unit current from
$a$ to $b$.
Traditionally this way of talking is reserved for time-reversible Markov chains,
which are precisely those for which we have the `reciprocity law'
$N_{abcd}=N_{cdab}$.
For such chains, if we build a resistor network
where nodes $i \neq j$ are
joined by a resistor of conductance (i.e., reciprocal resistance)
$-\Lap^{ij}$,
then $N_{abcd}$ will indeed
be the voltage drop between $c$ and $d$ due to a unit current from
$a$ to $b$.
We propose to extend this way of talking to non-time-reversible chains.

In electrical terms, the voltage drop $N_{abab}$ between $a$ and $b$ due to
a unit current between $a$ and $b$ is the \emph{effective resistance}.
This is the same as the reciprocal of the current that flows when a $1$-volt
battery is connected up between $a$ and $b$---which is what we get in effect
when we measure commuting times using green and red paint.
So the commuting time $C_{ab}=N_{abab}$ is the same as
the effective resistance between
$a$ and $b$.

The connection of commuting time to effective resistance,
and the general recognition
that commuting times play a key role in understanding Markov chains,
is due to Chandra et al.
\cite{crrst:commute}.

{\bf Note.}
Now we are in a position to understand the significance
of the name `cross-potential'.
This name is meant to indicate the connection of $N_{abcd}$
to the cross-ratio of complex function theory.
If we extend our notions about Markov chains to cover Brownian motion
on the Riemann sphere, we get
\begin{eqnarray*}
N_{abcd}
&=&
-\frac{1}{2\pi} (\log |a-c| - \log |a-d| - \log |b-c| + \log |b-d|)
\\&=&
-\frac{1}{2\pi} \log \left| \frac{a-c}{a-d} \frac{b-d}{b-c} \right|)
\\&=&
-\frac{1}{2\pi} \Re \log \frac{a-c}{a-d} \frac{b-d}{b-c}
.
\end{eqnarray*}
We don't have to specify a metric on the sphere here,
because the Laplacian is a conformal
invariant in two dimensions.
Thinking of the sphere as being an electrical conductor with constant
conductivity (say, 1 mho `per square'),
the electrical interpretation becomes exact.
The advantage of having $N$ to take four `arguments' now becomes apparent,
because 
$N_{abcb} = \infty$.
That's why engineers using look for cracks in nuclear reactor cooling pipes
with a \emph{4-point probe}.
To get a sensible generalization of $C_{ab}$ we will need to do some kind of
renormalization, which will introduce a dependence on the metric.
We should not be sorry about this,
because it brings curvature into the picture---and you know that can't be bad.

\section{Realization}

Now, finally, to realize commuting times as squared distances.
From the bilinear form $B$ we get the quadratic form
\[
Q(u)=||u||^2=B(u,u)
=
\sum_{ij} u^i Z_{ij} u^j
.
\]
\[
C_{ab}=N_{abab}=Q(\tdelta{_a^i}-\tdelta{_b}) = ||\tdelta{_a}-\tdelta{_b}||^2
.
\]
So if we map $i$ to $\tdelta{_i}$ then the commuting time $C_{ab}$
becomes the squared distance between the images in the $Q$-norm.

That is, if what we're calling the $Q$-norm is indeed a norm.
Is $Q$ really positive definite?

To understand better what is going on here, it is useful to look
at the bilinear form
\[
L(\phi,\psi)=
\sum_{ij} \phi_i \Lap^{ij} \psi_j
,
\]
where we think of $\phi$ and $\psi$ as being defined only modulo additive
constants.
If we think of $\phi_i$ as the potential of the measure
\[
\sum_i \phi_i \Lap^{ik} 
,
\]
then this is the same bilinear form as before, except that now instead
of measures of total mass $0$ it takes as its arguments the corresponding
potentials, the first with respect to the original chain, and the second
with respect to the time-reversed chain:
\[
L(\phi,\psi)
=
B(\sum_i \phi_i \Lap^{ik},
\sum_i \psi_i \Lap^{ki}
)
=
B(\sum_i \phi_i \Lap^{ik},
\sum_i \hat{\Lap}^{ik} \psi_i
)
.
\]
This follows from the formula $\Lap Z \Lap = \Lap$ above.

Now to get the equivalent of $Q$ in this context we restrict to the
subspace
\[
V=\{
(\phi,\psi): \sum_i \phi_i \Lap^{ik} = \sum_j \Lap^{kj} \psi_j
\}
\]
and take as our quadratic form
\[
R((\phi,\psi))=L(\phi,\psi)
.
\]

In the case of a time-reversible chain,
$V$ is just the diagonal $\phi = \psi$,
and
\[
Q(\phi \Lap)= R((\phi,\phi)) = L(\phi,\phi) =
\sum_{ij} \phi_i \Lap^{ij} \phi_j
= \half \sum_{ij} (-\Lap^{ij}) (\phi_i-\phi_j)^2
.
\]
This is evidently positive-definite.
Indeed, if we associate to $(\phi,\phi)$
the vector with $n \choose 2$ coordinates
$\sqrt{-\Lap^{ij}}(\phi_i-\phi_j)$, $i<j$,
then we will have embedded the normed space $(V,R)$,
and along with it our Markov chain,
in Euclidean $n \choose 2$-space.

Electrically, what we have done here is to account for the energy being
dissipated in the network by adding up the energy dissipated by
individual resistors.
And there should be some kind of probabilistic interpretation as well.

That's how it works for time-reversible chains, for which $\Lap^{ij}=\Lap^{ji}$.
However, the argument extends to the general case by what amounts to a trick.
The key is the observation that for $(\phi,\psi) \in V$ we have
\[
L(\phi,\psi)=L(\phi,\phi)=L(\psi,\psi)
.
\]
(But please note that in general $L(\phi,\psi) \neq L(\psi,\phi)$!)
So
\[
Q(\phi \Lap)= R((\phi,\psi)) =
L(\phi,\psi)=
L(\phi,\phi)=
\sum_{ij} \phi_i \Lap^{ij} \phi_j
= \half \sum_{ij} (-\Lap^{ij}) (\phi_i-\phi_j)^2
.
\]
So there is the positive-definiteness we need.

Now, though, we don't see any natural way to interpret
the terms of the sum electrically or probabilistically.
(Which is not to say that there isn't one!)
In putting $\phi$ in both slots of $L$ we leave the subspace $V$,
and thereby commit what appears to be an unnatural act.
But it seems to have paid off.

\section{Minimax characterization of commuting times and hitting probabilities}

Fix states $a \neq b$, and let
\[
S_{a,b}= \{ \phi | \phi_a =1, \phi_b=0 \}
\]
Here we really should be thinking of $\phi$ as being defined only up to an
additive constant, which means we should write $\phi_a-\phi_b=1$,
but we're going to be sloppy about this, because we want to focus attention
on two distinguished elements of $S_{a,b}$ which are naturally $1$ and $a$
and $0$ at $b$.
These are
\[
\phibar_i =
\Prob(\mbox{hit $a$ before $b$ starting at $i$ going backward in time})
\]
and
\[
\psibar_i =
\Prob(\mbox{hit $a$ before $b$ starting at $i$ going forward in time})
.
\]

We've met $\phibar$ before:  It's proportional to
the equilibrium concentration of green particles in our painting scenario.
$\psibar$ is the analogous quantity for the reversed chain.
The pair $(\phibar,\psibar)$ belongs to our subset $V$, because
\[
(\phibar \Lap)^i = (\Lap \psibar)^i =
r_{ab}(\tdelta{_a^i}-\tdelta{_b^i})
.
\]
Here we once again are writing $r_{ab} = \frac{1}{T_{ab}}$
for the equilibrium rate of commuting between $a$ and $b$.
Observe that any $f$ we have
\[
L(\phibar,f) = L(f,\psibar) = r_{ab} (f_a-f_b)
.
\]
So whenever $f$ is in $S_{a,b}$ we have
\[
L(\phibar,f) = L(f,\psibar) = r_{ab}
,
\]
and in particular
\[
L(\phibar,\psibar)=r_{ab}
.
\]

{\bf Theorem.}
\[
r_{ab}
=
\frac{1}{T_{ab}}
=
\min_{\alpha} \max_{\phi+\psi=2\alpha} L(\phi,\psi)
.
\]
Here and below,
$\alpha$, $\phi$, and $\psi$ are restricted to lie in $S_{a,b}$,
i.e. to take value $1$ at $a$ and $0$ at $b$.

{\bf Proof.}
Whatever $\alpha$ is,
we may take $\phi=\phibar$
(and thus $\psi=2\alpha-\phibar$),
and have
\[
L(\phi,\psi) = L(\phibar, \psi) = r_{ab}
\]
as above.
So
\[
\min_{\alpha} \max_{\phi+\psi=2\alpha} L(\phi,\psi)
\geq r_{ab}
.
\]

To prove the inequality in the other direction, and in the process
identify where the minimax is achieved,
take
\[
\alpha = \half(\phibar+\psibar)
.
\]
If $\phi+\psi=2\alpha$ then we can write
\[
\phi=\phibar+f
\]
and
\[
\psi=\psibar-f
,
\]
where $f_a=f_b=0$.

Now 
\[
L(\phibar,f)=L(f,\psibar)=r_{ab}(f_a-f_b)=0
,
\]
so
\[
L(\phi,\psi)
=
L(\phibar+f,\psibar-f)
=
L(\phibar,\psibar)-L(f,f)
= r_{ab}-L(f,f)
.
\]
And even though we claim it is a travesty to put the same $f$ into both
slots of $L$,
we still have
\[
L(f,f) \geq 0
:
\]
That was the upshot of our embedding investigation.
So
\[
L(\phi,\psi) \leq r_{ab}
,
\]
still assuming $\alpha=\half(\phibar+\psibar)$ and $\phi+\psi=2\alpha$.
Hence
\[
\min_{\alpha} \max_{\phi+\psi=2\alpha} L(\phi,\psi)
\geq r_{ab}
. \quad \qed
\]

In the time-reversible case, where $\Lap^{ij}=\Lap^{ji}$,
this minimax can be reduced to a straight
minimum.
That's because in this case for any $g,f$ we have
$L(f,g)=L(g,f)$, and hence
\[
L(g+f,g-f)
=
L(g,g)-L(f,f)
.
\]
So to maximize $L(\phi,\psi)$ while fixing the sum
$\phi+\psi=2\alpha$ we take $\phi=\psi=\alpha$.

{\bf Corollary.}
When $\Lap^{ij}$ is symmetric
\[
r_{ab} =
\min_{\phi(a)=1,\phi(b)=0} L(\phi,\phi)
.
\quad \qed
\]

This minimum principle for resistances was known already
to 19th century physicists, specifically Thomson (a.k.a. Kelvin),
Maxwell, and Rayleigh:
For more about this, see Doyle and Snell
\cite{doylesnell:walks}.

Having a straight minimum is a lot better than having a minimax,
because now we can plug in any $\phi$ with $\phi(a)=1,\phi(b)=0$
and get an upper bound for $r_{ab}$,
corresponding to a lower bound for $T_{ab}$.
This method is a staple of
electrical theory---the part of electrical theory that
doesn't extend to non-time-reversible chains because it depends on
the relation $L(f,g)=L(g,f)$.

For time-reversible chains there are also complementary methods for finding
lower bounds for $r_{ab}$, and thus upper bounds for $T_{ab}$.
These emerge from the minimum principle
through the mystery of convex duality.

\section{The obstruction to time-reversibility}

We close with an application to characterizing time-reversible chains.

Let $M_{ij}$ be the expected time to reach $j$ starting from $i$.
Coppersmith, Tetali, and Winkler
%\cite{ctw:obstruction}
showed that a Markov chain is time-reversible just if for all $a,b,c$
\[
M_{ab}+M_{bc}+M_{ca} =
M_{ac}+M_{cb}+M_{ba}
.
\]
And in this case the expected time to traverse a cycle of any length
will be the same in either direction.
Note that the $M_{ij}$s themselves are not conformally invariant,
these cycle sums are.
For a cycle of length $2$, the cycle sum is our best friend the commuting
time.

We always have
\[
M_{ab}+M_{bc}+M_{ca} =
\hat{M}_{ac}+\hat{M}_{cb}+\hat{M}_{ba}
\]
(look at a long record of the chain backwards),
so an equivalent condition is that for all $a,b,c$
\[
M_{ab}+M_{bc}+M_{ca} =
\hat{M}_{ab}+\hat{M}_{bc}+\hat{M}_{ca}
.
\]
This is true despite the fact that in general
\[
\hat{M}_{ab} \neq M_{ba}
.
\]

So, why is this true?  It comes down to the fact that a conformal class
of chains is reversible just if our bilinear form $L(\phi,\psi)$
on $V=\{x^i | \sum_i x_i = 0 \}$ is symmetric.
To any bilinear form $\sum_{ij} u^i Z_{ij} v^j$ on $V$ there corresponds
a natural cohomology class
\[
Z_{ij} - Z_{ji}
,
\]
which is to say, an antisymmetric matrix defined up to addition of a
matrix of the form $B_{ij}=a_i-a_j$.
This class represents the obstruction to symmetrizing the matrix of the form
within its $ab$-equivalence class.
This class vanishes just if it integrates to $0$ around any cycle,
and cycles of length $3$ span the space of cycles.
Indeed, they span it in a very redundant way.
To verify reversibility, it would suffice to check any basis for the space
of cycles,
e.g. only cycles of length $3$ involving the fixed state $n$
(the `ground').

\bibliography{commute}
\bibliographystyle{hplain}

\end{document}